\documentclass[11pt]{article}
\usepackage[T1]{fontenc}
\usepackage{a4wide}
\usepackage{lipsum}

\usepackage{lipsum}
\usepackage{amsmath,amssymb,amscd,amsthm,amsfonts,anyfontsize,color,enumerate,fix-cm,layout,lipsum,lpic,mathrsfs,mdwlist,stmaryrd,tensor,tikz,thmtools,xspace,comment}
\newtheorem{theorem}{Theorem}[section]

\theoremstyle{definition}           
\newtheorem{definition}[theorem]{Definition}

\newtheorem{problem}[theorem]{Problem}
\theoremstyle{remark}       
\newtheorem{remark}[theorem]{Remark} 

\usepackage{hyperref}

\title{Metric properties of electrical networks and the graph reconstruction problems}

\author{
V. G. Gorbounov \thanks{Faculty of Mathematics, National Research University Higher School of Economics, Usacheva 6, 119048 Moscow, Russia. E-mail: {\tt vgorb10@gmail.com}. }
\and
A.A. Kazakov \thanks{Lomonosov Moscow State University, Faculty of Mechanics and Mathematics, Russia, 119991, Moscow, GSP-1, 1 Leninskiye Gory, Main Building;
Kazan Federal University, N.I. Lobachevsky Institute of Mathematics and Mechanics,  Kazan, 420008, Russia. Supported by RSF (project No. 24-21-00158).}
}

\makeatletter

\newcommand{\shorttitle}{\@title}

\def\@maketitle{%
  \newpage
  \begin{center}%
  \let \footnote \thanks
    {\small Proceedings of the 13th European Conference on Combinatorics, Graph Theory and Applications\\ EUROCOMB'25\\
    Budapest, August 25 - 29, 2025
    }
    \vskip 0.5em
    \rule{\linewidth}{0.04cm}
    \vskip 3.5em
    {\LARGE \textbf{\textsc{\@title}} \par}%
    \vskip 1.5em
    {\textbf{\textsc{(Extended abstract)}} \par}
    \vskip 2.5em%
    {\large
      \lineskip .5em%
      \begin{tabular}[t]{c}%
        \@author
      \end{tabular}\par}%
  \end{center}%
  \par
  }
\makeatother

\usepackage{amsthm}

\begin{document}

\thispagestyle{empty}
\maketitle

\begin{abstract}
Using the generalized Temperley trick, we demonstrate the explicit embedding of circular electrical networks into totally non-negative Grassmannians. Building on this result, we show that the effective resistances between boundary nodes of circular electrical networks satisfy the Kalmanson property, and we provide the full characterization of planar electrical  Kalmanson metrics. Additionally, we present a graph reconstruction algorithm with applications in phylogenetic network analysis as well as the numerical solution of the Calderón problem.  
\\
\end{abstract}


\section{Introduction}

Electrical network theory is now a well-established applied field, originating in the mid-19th to early 20th centuries through the work of Kirchhoff, Maxwell, Kennelly, Norton, and many other physicists and engineers. From a mathematical  point of view, this theory is a great source of a large number of highly non-trivial combinatorial and geometric results and structures, which appear and find their applications in a wide variety of mathematical fields from   exactly solvable models  in statistical physics \cite{GoT}, \cite{Tal} to the famous  Calderon problem \cite{Kaz} and inverse problems in phylogenetic network theory \cite{GK}.

\section{Preliminaries}

\begin{definition} 
An \textbf{electrical network} $\mathcal{E}(G,  w)$ is a planar graph $G(V, E)$,  embedded into a disk and equipped with a  conductivity function $w: E(G) \to \mathbb{R}_{\geq 0}$, which together satisfy  the following conditions:  
\begin{itemize}
    \item All nodes are divided into the set of inner nodes $V_I$ and the set of boundary nodes $V_B$;
    \item An edge weight  $w(v_iv_j)=w_{ij}$  denotes the conductivity of this edge. 
\end{itemize}
An electrical network is called \textbf{circular} if   its boundary nodes lie on the boundary circle and are enumerated clockwise from $1$ to $|V_B|:=n$ (inner nodes are enumerated arbitrarily from $n+1$ to $|V|$).
\end{definition}

Consider an electrical network $\mathcal{E}(G, w)$ and apply voltages $\mathbf{U}: V_B \to \mathbb{R}$ to its boundary nodes $V_B$. These boundary voltages induce a unique harmonic extension $U: V \to \mathbb{R}$ to all vertices, which is determined by Ohm's and Kirchhoff's laws:
 \begin{equation*}
     \sum \limits_{j \in V}w_{ij}\bigr(U(i)-U(j)\bigl)=0, \ \forall i \in V_I.
 \end{equation*}

 One of the main objects associated with  each harmonic extension  of  boundary voltages   $\textbf{U}$ is the boundary currents $\textbf{I}=\{I_1, \dots, I_n\}$ running through boundary nodes:
 \begin{equation*}
    I_k:=\sum \limits_{j \in V}w_{ij}\bigr(U(k)-U(j)\bigl), \ k\in \{1, \dots, n\}.
 \end{equation*}

\begin{theorem} \textup{\cite{CIW}} \label{aboutresp}
Consider an electrical network $\mathcal{E}(G, w).$ Then, there is a  matrix $M_R(\mathcal{E})=(x_{ij}) \in Mat_{n \times n}(\mathbb{R})$ such that the following holds:
$$M_R(\mathcal{E})\textbf{U}=\textbf{I}.$$
This matrix is called the \textbf{response matrix} of a network $\mathcal{E}(G, \omega).$ 
\end{theorem} 
Taking the <<inverse>> of a matrix $M_R(\mathcal{E})$, we can obtain another important electrical network characteristic:
\begin{definition} \label{lem:eff-resist}
Let $\mathcal E(G, \omega)$ be a connected  electrical network with $n$ boundary nodes, and let the boundary voltages  $U = (U_1, \dots , U_n)$ be such that
\begin{equation*} \label{eq-resist}
    M_R(\mathcal E)U = -e_i + e_j,
\end{equation*}
 where $e_k, \ k \in \{1, \dots  , n\}$ is the  standard basis of $\mathbb{R}^n$. 
 
 Let us define the \textbf{effective resistance} $R_{ij}$ between  nodes $i$ and $j$ as follows $|U_i - U_j|:=R_{ij}.$   Effective resistances are well-defined and $R_{ij}=R_{ji}$. 

We will organize the effective resistances $R_{ij}$ in an  \textbf{effective resistances matrix}  $R(\mathcal E)$ setting $R_{ii}=0$ for all $i$.
\end{definition} 
A well-known result, with broad applications in applied mathematics and organic chemistry, states that:
\begin{theorem} \textup{\cite{Kl}} \label{th:about metric}
    Let $\mathcal E(G, \omega)$ be a connected  electrical network  then for any of its three boundary nodes $k_1, k_2$ and $k_3$ the triangle inequality holds:
    $$R_{k_1k_3}+R_{k_2k_3}-R_{k_1k_2} \geq 0.$$
    Hence the set of all $R_{ij}$ defines a metric on the boundary nodes of $G$.
\end{theorem}

\section{Main results}
Kirchhoff's classical matrix tree theorem helps us to calculate $R_{ij}$ using partition functions of spanning trees. Building upon this, the generalized Temperley trick \cite{L} allows us to alternatively calculate $R_{ij}$ with dimer (almost perfect matching) partition functions. The last combinatorial statement is the key to explicitly relating electrical network theory with the geometry of totally non-negative Grassmannians  $Gr_{\geq 0}(n-1,2n)$:

    \begin{theorem}  \textup{\cite{BG1}, \cite{GK}} \label{main-gr}
Let $\mathcal E(G, \omega)$ be a connected electrical network with $n$ boundary nodes. Using its effective resistances matrix $R(\mathcal E)$,  let us define a  matrix:
\begin{equation} \label{eq:omega_n,r}
  \Omega_{R}(\mathcal E)=\left(\begin{matrix}
1 & m_{11} & 1 &  -m_{12} & 0 & m_{13} & 0 & \ldots  \\
0 & -m_{21} & 1 & m_{22} & 1 & -m_{23} & 0 & \ldots \\
0 & m_{31} & 0 & -m_{32} & 1 & m_{33} & 1 & \ldots \\
\vdots & \vdots & \vdots & \vdots & \vdots & \vdots &  \vdots & \ddots 
\end{matrix}\right),
\end{equation}
where 
$$m_{ij}= -\frac{1}{2}(R_{i,j}+R_{i+1,j+1}-R_{i,j+1}-R_{i+1,j}).$$ 

Then, $\Omega_{R}(\mathcal E)$ defines a  point $\mathcal{L}(\mathcal{E})$ in the totally non-negative part of  $Gr(n-1,2n)$, it means that:
\begin{itemize}
    \item  The dimension of the row space of $\Omega_R(\mathcal{E})$ is equal to $n-1$;
    \item Each $n-1 \times n-1$ minors     of  $\Omega_R(\mathcal{E})$  is non-negative;
   \item   Plucker coordinates of the point of $\mathcal{L}(\mathcal{E})$ correspond  to $n-1 \times n-1$ minors  of the matrix $\Omega_R'(\mathcal{E})$ obtained from $\Omega_R(\mathcal{E})$  by deleting the last row.
\end{itemize}
\end{theorem}
Many interesting inequalities involving $R_{ij}$ follow from the positivity of the Plucker coordinates of the point represented by $\Omega_{R}(\mathcal E)$. Some of them have an elegant interpretation:

\begin{theorem}  \textup{\cite{GK}} \label{charkalm}
     Let $\mathcal E(G, \omega)$ be a connected circular electrical network and let  $i_1, i_2, i_3, i_4$ be any four nodes in the circular order. Then the   Kalmanson inequalities hold:
      \begin{equation*} \label{kal_1}
      R_{i_1i_3}+R_{i_2i_4}\geq \textit{max}(R_{i_2i_3}+R_{i_1i_4},  R_{i_1i_2}+R_{i_3i_4}).   
      \end{equation*}

    Since,  $R_{ij}$ can be considered as the Kalmanson metric on boundary nodes. 
\end{theorem}
Additionally, Theorem \ref{main-gr} allows us to obtain the exhaustive  characterization of planar electrical  Kalmanson metrics:
\begin{theorem} \textup{\cite{GK}} \label{th: dual}
  Let  $D=(d_{ij}) \in \text{Mat}_{n \times n}(\mathbb{R})$ be a matrix of a Kalmanson metric with respect of a cyclic numeration $1, \dots, n$, then $D$ is the effective resistance matrix of a connected circular electrical network $\mathcal E$ with $n$ boundary nodes if and only if  the matrix $\Omega_{D}$ constructed from $D$ according to the formula \eqref{eq:omega_n,r} defines a point $X$ in $\mathrm{Gr}_{\geq 0}(n-1, 2n)$ and the Plucker coordinate $\Delta_{24\dots 2n-2}(X)$ does not vanish.
\end{theorem}
The last theorem can be elegantly reformulated in the following equivalent form:
\begin{theorem} \textup{\cite{GK}} \label{th: dual1}
  Let  $D=(d_{ij}) \in \text{Mat}_{n \times n}(\mathbb{R})$ be a matrix of a Kalmanson metric with respect to a cyclic numeration $1, \dots, n$ and let us consider a matrix $M(D)=(m_{ij}),$ $m_{ij}=\frac{1}{2}(d_{i,j}+d_{i+1,j+1}-d_{i,j+1}-d_{i+1,j}) $ of the coefficients in its   circular split decomposition, see \cite{DL}. Then  $D$  is the effective resistance matrix of a connected circular electrical network $\mathcal E$ with $n$ boundary nodes if and only if  the matrix    $-M(D)$ is a response matrix $M_R(\mathcal{E^{*}})$ of a network $\mathcal{E^{*}}$. Moreover, $\mathcal{E^{*}}$  is the  dual network to $\mathcal{E}.$
\end{theorem}

\section{Applications}
A notoriously difficult problem in electrical network theory is to recover a graph of an electrical network from its effective resistance matrix, see \cite{CIW}. Beyond its fundamental applications in physics and engineering (see \cite{BDV}, \cite{Kaz}, \cite{Uh2}), solutions to this problem can also be used for the reconstruction of phylogenetic networks (see \cite{F1}, \cite{F2}, \cite{GK}), owing to the fact that the effective resistance metric of circular networks satisfies the Kalmanson property. We will focus on presenting a sketch of a new algorithm for reconstructing the topology of a graph of a circular electrical network from its effective resistance matrix.
\begin{problem}  \label{bl-box-res-top}
    Consider an electrical network $\mathcal{E}(G, w)$ on an unknown  graph $G$.  It is required to reconstruct $($up to electrical transformations, see \cite{CIW}$)$ a graph $G$   by a known effective resistance matrix $R_{\mathcal E}$. 
\end{problem}

The ability to solve the last problem is revealed by the following results.
\begin{definition} \textup{\cite{K}}
A circular electrical network is called \textbf{minimal} if the strands of its median graph do not have self-intersections; any two strands intersect at most one point, and the median graph has no loops or lenses.
\end{definition}
\begin{figure}[h!]
    \centering
    \includegraphics[width=0.4\textwidth]{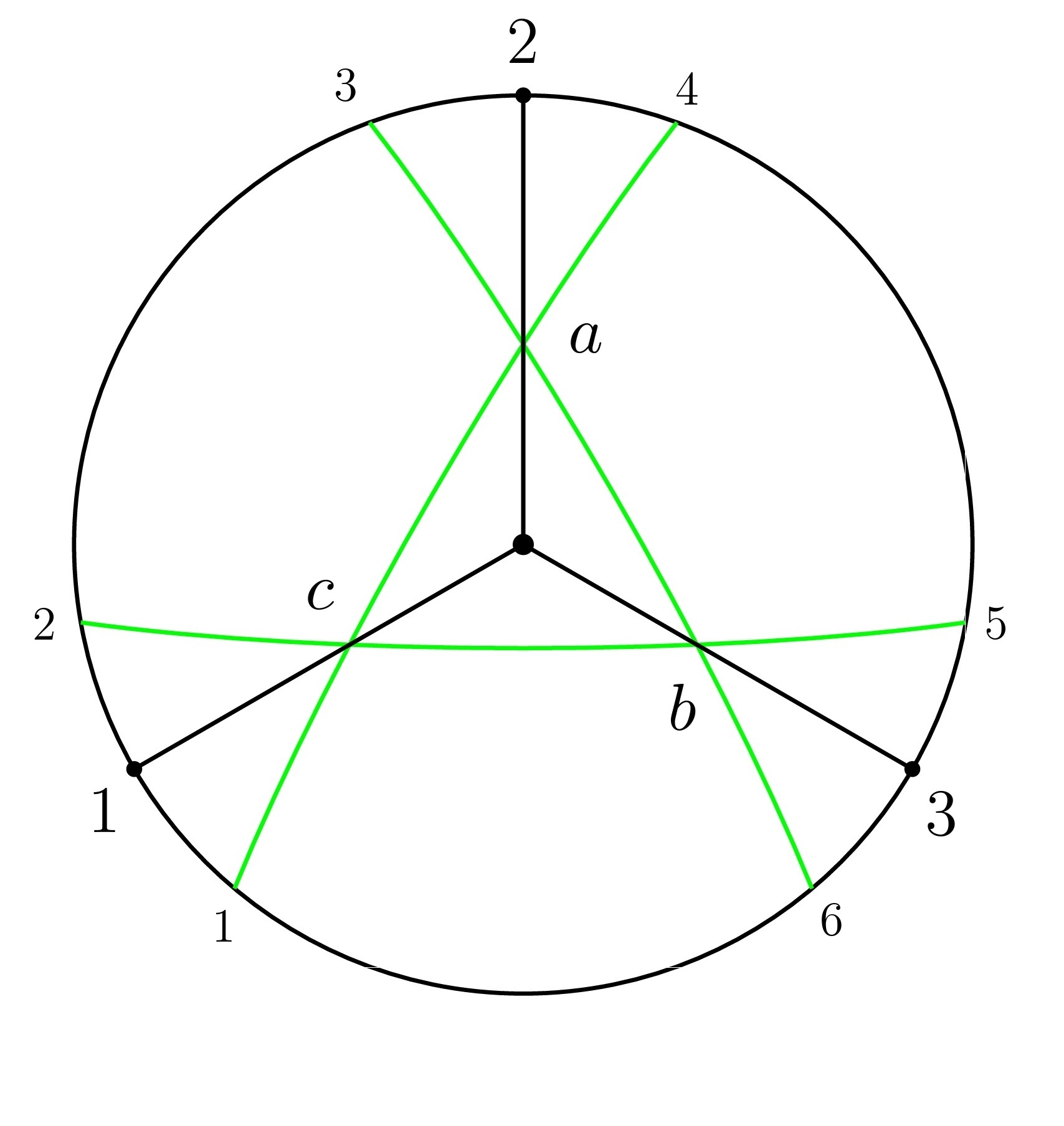}
    \caption{Star-shape network, its median graph and the strand permutation $\tau(\mathcal E)=(14)(36)(25)$ }
    \label{fig:triangle}
\end{figure}
\begin{theorem} \textup{\cite{CIW}}
    Each minimal circular electrical network $\mathcal E(G, \omega)$  is unique up to electrical transformations defined by its strand permutations $\tau(\mathcal{E}).$
\end{theorem}
Denote by $A_i$ the columns of the matrix $\Omega_{ R}(\mathcal E)$ and define the column permutation $g(\mathcal E)$ as follows:
    $g(\mathcal E)(i)=j,$ if  $j$  is the minimal number such that   $A_i \in \ \mathrm{span}(A_{i+1}, \dots, A_{j}  ),$ where the indices are taken modulo $2n$.

\begin{theorem} \textup{\cite{GK}} \label{bl-box_th-res}	
Up to the star-triangle transformations, a topology of each  minimal  electrical network $\mathcal{E}(G, w)$ with $n$ boundary nodes  can be uniquely recovered by a column  <<rank-patterns>>  of $\Omega_R(\mathcal{E})$:
$$g(\mathcal E)+1=\tau(\mathcal E) \mod 2n.$$
\end{theorem}

\begin{remark}
    A more advanced technique called the generalized chamber ansatz  applied to $\Omega_{ R}(\mathcal E)$ provides an algorithm for recovering  the conductivity function $\omega$  as well, see \cite{Kaz}.  
\end{remark}


\end{document}